\documentclass[12pt]{article}
\usepackage{graphicx}
\usepackage{amsmath,amsthm,amssymb,enumerate}
\usepackage{euscript,mathrsfs}
\usepackage[left=2cm,right=2cm,top=3.5cm,bottom=3.5cm]{geometry}
\usepackage{color}

\usepackage{soul}

\catcode`\@=11 \@addtoreset{equation}{section}

\catcode`\@=12

\allowdisplaybreaks

\newtheorem{Theorem}{Theorem}[section]
\newtheorem{Proposition}[Theorem]{Proposition}
\newtheorem{Lemma}[Theorem]{Lemma}
\newtheorem{Corollary}[Theorem]{Corollary}

\theoremstyle{definition}
\newtheorem{Definition}[Theorem]{Definition}

\newtheorem{Remark}[Theorem]{Remark}

\newcommand{\bTheorem}[1]{
\begin{Theorem} \label{T#1} }
\newcommand{\eT}{\end{Theorem}}

\newcommand{\bProposition}[1]{
\begin{Proposition} \label{P#1}}
\newcommand{\eP}{\end{Proposition}}

\newcommand{\bLemma}[1]{
\begin{Lemma} \label{L#1} }
\newcommand{\eL}{\end{Lemma}}

\newcommand{\bCorollary}[1]{
\begin{Corollary} \label{C#1} }
\newcommand{\eC}{\end{Corollary}}

\newcommand{\bRemark}[1]{
\begin{Remark} \label{R#1} }
\newcommand{\eR}{\end{Remark}}

\newcommand{\Td}{\mathbb{T}^d}

\newcommand{\bDefinition}[1]{
\begin{Definition} \label{D#1} }
\newcommand{\eD}{\end{Definition}}

\newcommand{\vrd}{\varrho_\delta}
\newcommand{\vud}{\vu_\delta}
\newcommand{\chid}{\chi_\delta}

\newcommand{\Del}{\Delta_x}

\newcommand{\Ds}{\mathbb{D}_x}

\newcommand{\bfphi}{\boldsymbol{\varphi}}

\newcommand{\bFormula}[1]{
\begin{equation} \label{#1}}
\newcommand{\eF}{\end{equation}}

\newcommand{\Ov}[1]{\overline{#1}}

\newcommand{\DC}{C^\infty_c}

\newcommand{\aleq}{\stackrel{<}{\sim}}

\newcommand{\ageq}{\stackrel{>}{\sim}}

\newcommand{\vr}{\varrho}
\newcommand{\vre}{\vr_\ep}

\newcommand{\vue}{\vu_\ep}

\newcommand{\vu}{\vc{u}}
\newcommand{\vm}{\vc{m}}

\newcommand{\vc}[1]{{\bf #1}}

\newcommand{\Div}{{\rm div}_x}
\newcommand{\Grad}{\nabla_x}

\newcommand{\dx}{\,{\rm d} {x}}

\newcommand{\dt}{\,{\rm d} t }

\newcommand{\intTd}[1]{\int_{\Td} #1 \ \dx}
\newcommand{\intO}[1]{\int_{\Omega} #1 \ \dx}

\newcommand{\D}{{\rm d}}

\newcommand{\ep}{\varepsilon}

\def\softd{{\leavevmode\setbox1=\hbox{d}%
          \hbox to 1.05\wd1{d\kern-0.4ex{\char039}\hss}}}
\definecolor{Cgrey}{rgb}{0.85,0.85,0.85}
\definecolor{Cblue}{rgb}{0.50,0.85,0.85}
\definecolor{Cred}{rgb}{1,0,0}
\definecolor{fancy}{rgb}{0.10,0.85,0.10}

\newcommand\Cbox[2]{%
    \newbox\contentbox%
    \newbox\bkgdbox%
    \setbox\contentbox\hbox to \hsize{%
        \vtop{
            \kern\columnsep
            \hbox to \hsize{%
                \kern\columnsep%
                \advance\hsize by -2\columnsep%
                \setlength{\textwidth}{\hsize}%
                \vbox{
                    \parskip=\baselineskip
                    \parindent=0bp
                    #2
                }%
                \kern\columnsep%
            }%
            \kern\columnsep%
        }%
    }%
    \setbox\bkgdbox\vbox{
        \color{#1}
        \hrule width  \wd\contentbox %
               height \ht\contentbox %
               depth  \dp\contentbox
        \color{black}
    }%
    \wd\bkgdbox=0bp%
    \vbox{\hbox to \hsize{\box\bkgdbox\box\contentbox}}%
    \vskip\baselineskip%
}

\date{}

\begin{document}


\title{Two phase flows of compressible viscous fluids}

\author{Eduard Feireisl
\thanks{The work of E.F. was supported by the
Czech Sciences Foundation (GA\v CR), Grant Agreement
21--02411S.}
\and Anton\' \i n Novotn\' y { \thanks{The work of A.N. was partially supported by the Eduard \v Cech visiting program at the Mathematical Institute of the Acdemy of Sciences of the Czech Republic. }
}
}


\maketitle

\centerline{Institute of Mathematics of the Academy of Sciences of the Czech Republic;}
\centerline{\v Zitn\' a 25, CZ-115 67 Praha 1, Czech Republic}

\centerline{feireisl@math.cas.cz}

\centerline{and}

\centerline{IMATH, EA 2134, Universit\'e de Toulon,}
\centerline{BP 20132, 83957 La Garde, France}
\centerline{novotny@univ-tln.fr}

\begin{abstract}
	
	We introduce a new concept of \emph{dissipative varifold solution} to models of two phase compressible viscous fluids. In contrast with the existing approach based on the Young measure description, the new formulation is variational combining the energy and momentum balance in a single inequality. We show the existence of dissipative varifold solutions for a large class of general viscous fluids with non--linear dependence of the viscous stress on the symmetric velocity gradient.

\end{abstract}


\bigskip

\centerline{\it Dedicated to Maurizio Grasselli on the occasion of his 60-th birthday}

\bigskip


\section{Introduction}
\label{p}

We consider a simple model of a two--phase flow, where the interface between the two fluids is described via a phase variable $\chi \in \{ 0 , 1 \}$ that coincides with the characteristic function of the domain occupied by one of the fluid 
components. Accordingly, the time evolution of $\chi$ is governed by the transport equation 
\begin{equation} \label{p0}
	\partial_t \chi + \vu \cdot \Grad \chi = 0,
	\end{equation}
where $\vu$ is the fluid velocity.
We denote 
\[
\begin{split}
\mathcal{O}_1(t) &= \left\{  x \in \Omega \ \Big| \ \chi(t,x) = 1 \right\} - \mbox{the part of the physical space}\ 
\Omega\subset R^d \ \mbox{occupied by fluid 1},\\
\mathcal{O}_2(t) &= \left\{  x \in \Omega \ \Big| \ \chi(t,x) = 0 \right\} - \mbox{the part of the physical space}\ 
\Omega\subset R^d \ \mbox{occupied by fluid 2},\\ 
\Gamma(t) &= \Ov{\mathcal{O}}_1(t) \cap \Ov{\mathcal{O}}_2 (t) - \mbox{the interface.} 
\end{split}
\]
As the mass of the fluid is conserved, its density $\vr$ satisfies the standard equation of continuity
\begin{equation} \label{p1}
\partial_t \vr + \Div (\vr \vu ) = 0.
\end{equation}

The material properties of the two fluid components are characterized by a general barotropic equation of state  
\begin{equation} \label{p2a}
p_i = p_i(\vr),\ i=1,2,\ p(\chi, \vr) = \chi p_1 (\vr) + (1 - \chi)p_2(\vr). 
\end{equation}

The viscous stress $\mathbb{S}$ is related to the symmetric velocity gradient $\Ds \vu$ through an ``implicit'' rheological law
\begin{equation} \label{p2}
\mathbb{S}_i : \Ds \vu = F_i (\Ds \vu) + F^*_i (\mathbb{S}_i),\ i = 1,2	
\end{equation}
where the dissipative potentials $F_i$, $i=1,2$ are convex functions defined on the space $R^{d \times d}_{\rm sym}$ of 
symmetric $d-$dimensional real matrices, and where $F^*_i$ denotes the convex conjugate.
More specifically, we suppose 
\begin{equation} \label{p4}
	F_i : R^{d \times d}_{\rm sym} \to [0, \infty] \ \mbox{is l.s.c. convex},\ F_i(0) = 0,\ i=1,2.
\end{equation}
Note that \eqref{p2} is equivalent  
to
\[
\mathbb{S}_i \in \partial F_i (\Ds \vu) \ \Leftrightarrow \ \Ds \vu \in \partial F^*_i (\mathbb{S}_i).
\]
It is easy to check that the choice 
\[
F_i = \frac{\mu_i}{2} \left| \mathbb{D} - \frac{1}{d} {\rm tr} [\mathbb{D}] \right|^2 + 
\frac{\lambda_i}{2} \left|{\rm tr} \mathbb[\mathbb{D}] \right|^2, \ i=1,2
\]
gives rise to the standard Newtonian viscous stress used in the Navier--Stokes system. 
Similarly to \eqref{p2a}, we set 
\[
\mathbb{S} (\chi, \Ds \vu) = \chi \mathbb{S}_1 (\Ds \vu) + (1 - \chi) \mathbb{S}_2 (\Ds \vu).
\]
We refer to Bul\' \i \v cek et al 
\cite{BuGwMaSG} for more details about ''implicitly'' constituted fluids, where the viscous stress is related to the velocity gradient in a way similar to \eqref{p2}. Similar abstract approach has been used in \cite{AbbFeiNov}.

The balance of momentum is satisfied for each individual fluid. The velocity is continuous, while the Cauchy stress
\[
\mathbb{T}_i = \mathbb{S}_i - p_i \mathbb{I},\ i=1,2
\]
experiences a jump on the interface,
specifically,
\[
[\mathbb{T} \cdot \vc{n}] = \kappa H \vc{n}\ \mbox{on} \ \Gamma ,\ \kappa \geq 0,
\]
where $H \vc{n}$ is the mean curvature vector, and $\kappa$ the
coefficient of surface tension.
Consequently, it is convenient to formulate the momentum balance in the weak form:
\begin{equation} \label{p1a}
	\begin{split}
\left[ \intO{\vr \vu \cdot \bfphi } \right]_{t = 0}^{t = \tau} &= \int_0^\tau \intO{ \Big( \vr \vu \cdot \partial_t \bfphi + \vr \vu \otimes \vu: \Grad \bfphi + p(\chi, \vr) \Div \bfphi \Big) } \dt \\ &- \int_0^\tau \intO{ \mathbb{S} (\chi , \Grad \vu) : \Grad \bfphi } \dt 
+ \kappa \int_0^\tau \int_{\Gamma(t)} H \vc{n} \cdot \bfphi \ \D S_x \dt 
\end{split}
\end{equation}
for any $0 \leq \tau \leq T$, $\bfphi \in C^1_c([0,T] \times \Omega; R^d)$,
where we have omitted the effect of external volume forces for simplicity. 

\subsection{Varifold solutions}

If the initial interface is a smooth surface (curve if $d=2$), the problem admits a local in time regular solution as long as the data are regular, see the survey of Denisova, Solonnikov \cite{DeSo} and the references cited therein. Global existence is not 
expected because of possible singularities and self--intersections of $\Gamma (t)$.  

Plotnikov \cite{Plot1, Plot3, Plot2} adapted the concept of \emph{varifold} to study the problem for large times in the context of 
incompressible non-Newtonian fluids for $d=2$. Later Abels \cite{AbelsI, AbelsII} introduced the class of \emph{measure--valued} 
varifold solutions to attack the problem for $d=2,3$ still in the incompressible setting but for a larger family of 
viscous stresses including the standard Newtonian fluids. He also showed the existence of weak varifold solutions in the absence of surface tension. Ambrose at al. \cite{ALNS} extended the existence proof to the case with surface tension for incompressible Newtonian fluids and $d=2$. To the best of our knowledge, 
with the only exception of the work of Plotnikov \cite{Plot2} dealing with a simplified Stokes model, similar problems are completely open in the context of \emph{compressible} fluids.
	
We introduce a new weak formulation of the problem based on the principles of the calculus of variations. The main idea is to rewrite the momentum equation with the associated total energy balance as a single inequality in terms of the dissipation potentials $F_1$, $F_2$. The resulting problem still uses the concept of varifold, however, the description of oscillations via a Young measure is no longer necessary. The new formulation is versatile and can be used for both compressible and incompressible fluids. Finally, we show that the problem admits global in time solutions as long as $\Div \vu$ is penalized in a way similar to \cite{FeLiMa}.

\section{Weak formulation}

To avoid problems with physical boundary, we consider the space--periodic boundary conditions. In other words, 
the spatial domain $\Omega$ is identified with a flat torus, 
\begin{equation} \label{p5a}
	\Omega = \Td,\ d=2,3.
	\end{equation}

\subsection{Varifolds}

Following Abels \cite{AbelsI}, Plotnikov \cite{Plot1}, we introduce the concept of \emph{varifold} $V$, 
\[
V \in L^\infty_{\rm weak-(*)}(0,T; \mathcal{M}^+(\Td \times S^{d-1})),
\]
which, after disintegration, takes the form 
\[
V = |V| \otimes V_x,\ |V| \in \mathcal{M}^+(\Td),\ \{ V_x \}_{x \in \Td} \subset \mathfrak{P}(S^{d-1}).
\]
The \emph{first variation of} $V$ reads 
\[
\left< \delta V(t); \bfphi \right> = \int_{\Td \otimes S^{d-1}} (\mathbb{I} - z \otimes z) : \Grad \bfphi(x) \ \D V(t)
\ \mbox{for any}\ \bfphi \in C^1(\Td).
\]
Similarly to \cite{AbelsI}, we rewrite the momentum equation \eqref{p1a} as
\begin{equation} \label{p1aa}
	\begin{split}
		\left[ \intTd{\vr \vu \cdot \bfphi } \right]_{t = 0}^{t = \tau} &= \int_0^\tau \intTd{ \Big( \vr \vu \cdot \partial_t \bfphi + \vr \vu \otimes \vu: \Grad \bfphi + p(\chi, \vr) \Div \bfphi \Big) } \dt \\ &- \int_0^\tau \intTd{ \mathbb{S} (\chi , \Ds \vu) : \Grad \bfphi } \dt 
		+ \kappa \int_0^\tau \int_{\Td \times S^{d-1}}  \left< \delta V; \bfphi \right>  \ \D V \dt 
	\end{split}
\end{equation}
for any $0 \leq \tau \leq T$, $\bfphi \in C^1([0,T] \times \Td; R^d)$.

\subsection{Total energy balance}

Introducing the pressure potentials
\begin{equation} \label{p5b}
P_i(\vr),\ P'_i(\vr) \vr - P_i(\vr) = p_i(\vr),\ i=1,2,\ P(\chi,\vr) = \chi P_1(\vr) + (1-\chi) P_2(\vr),
\end{equation}	
we deduce the \emph{total energy balance} 
\begin{equation} \label{p5c}
\frac{{\rm d}}{\dt } \left[ \intTd{ \left( \frac{1}{2} \vr|\vu|^2 + P(\chi, \vr) \right) } + \kappa \| \Grad \chi \|_{\mathcal{M}(\Td; R^d)} \right] + \intTd{ \mathbb{S}(\chi, \Ds \vu) : \Ds \vu } = 0,
\end{equation}
cf. Abels \cite{AbelsI}.
In the weak formulation, the energy balance \eqref{p5c} is replaced by an inequality 
\begin{equation} \label{p5d}
\left[ \intTd{ \left( \frac{1}{2} \vr|\vu|^2 + P(\chi, \vr) \right) } + \kappa \| \Grad \chi \|_{\mathcal{M}(\Td; R^d)} \right]_{t=0}^{t=\tau} + \int_0^\tau \intTd{ \mathbb{S}(\chi, \Ds \vu) : \Ds \vu } \dt \leq 0.
	\end{equation}

Next, subtracting \eqref{p1aa} from \eqref{p5d} we obtain 
\begin{equation} \label{p5e}
	\begin{split}
	&\left[ \intTd{ \left( \frac{1}{2} \vr|\vu|^2  - \intTd{\vr \vu \cdot \bfphi } + P(\chi, \vr) \right) } + \kappa \| \Grad \chi \|_{\mathcal{M}(\Td; R^d)} \right]_{t=0}^{t=\tau}\\
	 &+ \int_0^\tau \intTd{ \Big( \vr \vu \cdot \partial_t \bfphi + \vr \vu \otimes \vu: \Grad \bfphi + p(\chi, \vr) \Div \bfphi \Big) } \dt \\ & 
	+ \kappa \int_0^\tau \int_{\Td \times S^{d-1}}  \left< \delta V; \bfphi \right>  \ \D V \dt \\
	&+ \int_0^\tau \intTd{ \mathbb{S}(\chi, \Ds \vu) : (\Ds \vu - \Ds \bfphi) } \dt \leq 0.
	\end{split}
\end{equation}

Finally, in view of \eqref{p2},
\[
\begin{split}
F(\chi, \Ds \bfphi) - F(\chi, \Ds \vu) &\equiv \chi \Big( F_1(\Ds \bfphi ) - F_1 (\Ds \vu) \Big) + (1 -\chi) \Big( F_2(\Ds \bfphi ) - F_2 (\Ds \vu) \Big)\\ 
&\geq \mathbb{S}(\chi, \Ds \vu) : (\Ds \bfphi - \Ds \vu) 
\end{split}
\]
Consequently, we rewrite \eqref{p5e} in the final form 
\begin{equation} \label{p5f}
	\begin{split}
		\int_0^\tau &\intTd{ \Big( F(\chi, \Ds \bfphi) - F(\chi, \Ds \vu) \Big) }  \dt\\ &\geq
		\left[ \intTd{ \left( \frac{1}{2} \vr|\vu|^2  - \intTd{\vr \vu \cdot \bfphi } + P(\chi, \vr) \right) } + \kappa \| \Grad \chi \|_{\mathcal{M}(\Td; R^d)} \right]_{t=0}^{t=\tau}\\
		&+ \int_0^\tau \intTd{ \Big( \vr \vu \cdot \partial_t \bfphi + \vr \vu \otimes \vu: \Grad \bfphi + p(\chi, \vr) \Div \bfphi \Big) } \dt \\ & 
		+ \kappa \int_0^\tau \int_{\Td \times S^{d-1}}  \left< \delta V; \bfphi \right>  \ \D V \dt
	\end{split}
\end{equation}
for any $0 \leq \tau \leq T$, $\bfphi \in C^1([0,T] \times \Td; R^d)$.

\subsection{Dissipative varifold solutions}

We are ready to introduce the concept of \emph{dissipative varifold solution}.

\begin{Definition}[{\bf dissipative varifold solution}] \label{pD1}
	We say that $[\chi, \vr, \vu, V]$ is a \emph{dissipative varifold solution} of the problem
	\eqref{p0}, \eqref{p1}, \eqref{p1a}, \eqref{p5a} if the following holds: 
	\begin{itemize}
		\item {\bf Regularity class.}
		\begin{equation} \label{ds1}
			\chi \in C([0,T]; L^1(\Td)),\ \chi \in \{ 0,1 \} \ \mbox{for a.a.}\ t, x;
			\end{equation}
		\begin{equation} \label{ds2}
			\vr \in C([0,T; L^1(\Td)]),\ 0 \leq \underline{\vr} \le \vr(t,x) \leq \Ov{\vr}\ 
			\mbox{for a.a.}\ t, x;
		\end{equation}
	\begin{equation} \label{ds3}
		\vu \in L^\alpha (0,T; W^{1,\alpha} (\Td; R^d))\ \mbox{for some}\ \alpha > 1, 
		\vr \vu \in C([0,T]; L^2(\Td; R^d)),
		\end{equation}
	\begin{equation} \label{ds4}
		V \in L^\infty_{\rm weak-(*)}(0,T; \mathcal{M}^+ (\Td \times S^{d-1})).
	\end{equation}

\item {\bf Transport.}
\begin{equation} \label{ds5}
	\left[ \intTd{ \chi \varphi } \right]_{t=0}^{t = \tau} = 
	\int_0^\tau \intTd{ \Big( \chi \partial_t \varphi + \chi \vu \cdot \Grad \varphi + \chi \Div \vu \varphi \Big) } \dt
	\end{equation}
	 for any $0 \leq \tau \leq T$, $\varphi \in C^1([0,T \times \Td])$. 
	 \item {\bf Mass conservation.}
	 \begin{equation} \label{ds6}
	 	\left[ \intTd{ \vr \varphi } \right]_{t=0}^{t = \tau} = 
	 	\int_0^\tau \intTd{ \Big( \vr \partial_t \varphi + \vr \vu \cdot \Grad \varphi \Big) } \dt
	 \end{equation}
 for any $0 \leq \tau \leq T$, $\varphi \in C^1([0,T \times \Td])$.
 \item{{\bf Momentum--energy balance}}
 \begin{equation} \label{ds7}
 	\begin{split}
 		\int_0^\tau &\intTd{ \Big( F(\chi, \Ds \bfphi) - F(\chi, \Ds \vu) \Big) }  \dt\\ &\geq
 		\left[ \intTd{ \left( \frac{1}{2} \vr|\vu|^2  - \vr \vu \cdot \bfphi  + P(\chi, \vr) \right) } + \kappa \| \Grad \chi \|_{\mathcal{M}(\Td; R^d)} \right]_{t=0}^{t=\tau}\\
 		&+ \int_0^\tau \intTd{ \Big( \vr \vu \cdot \partial_t \bfphi + \vr \vu \otimes \vu: \Grad \bfphi + p(\chi, \vr) \Div \bfphi \Big) } \dt \\ & 
 		+ \kappa \int_0^\tau \int_{\Td \times S^{d-1}}  \left< \delta V; \bfphi \right>  \ \D V \dt
 	\end{split}
 \end{equation}
 for any $0 \leq \tau \leq T$, and any $\bfphi \in C^1([0,T] \times \Td; R^d)$ satisfying 
 \begin{equation} \label{ds7a}
 \int_0^T \intTd{ \left( F_1 (\Ds \bfphi) + F_2 (\Ds \bfphi) \right) } \dt < \infty.
 \end{equation}
 \item {\bf Varifold compatibility.}
 If, in addition,  
 $\kappa > 0$, we require $\chi \in L^\infty_{{\rm weak-(*)}}(0,T; BV(\Td))$, and
 \begin{equation} \label{ds8}
 \int_{\Td \times S^{d-1}} \bfphi \cdot z \ \D V(\tau, \cdot) + \int_{\Td} \bfphi \cdot \D \Grad \chi (\tau) = 0	
 	\end{equation}
 for a.a. $0 \leq \tau \leq T$, and any $\bfphi \in C(\Td)$.
 
 \end{itemize}
	\end{Definition}

First observe that dissipative varifold solutions coincide with weak varifold solutions in the sense of Plotnikov \cite{Plot1} and Abels \cite{AbelsI}, if they satisfy the energy balance 
\begin{equation} \label{ds9}
\left[ \intTd{ \left( \frac{1}{2} \vr |\vu|^2 + P(\chi, \vr) \right) } + 
\kappa \| \Grad \chi \|_{\mathcal{M}(\Td; R^d)} \right]_{t=0}^{t = \tau} + 
\int_0^\tau \intTd{ \mathbb{S} : \Ds \vu } \dt = 0,  	 
	\end{equation} 
	where 
\begin{equation} \label{ds10}
	\mathbb{S}(t,x) \in \chi \partial F_1((\chi , \Ds \vu) (t,x)) + (1- \chi) \partial F_2((\chi , \Ds \vu) (t,x)) \ \mbox{for a.a.}\ t,x.
	\end{equation}
Indeed subtracting \eqref{ds9} from \eqref{ds7} we obtain 
\[
	\begin{split}
		\int_0^\tau &\intTd{ \Big( F(\chi, \Ds \bfphi) - F(\chi, \Ds \vu) - \mathbb{S} : (\Ds \bfphi - \Ds \vu \Big) }  \dt\\ &\geq
		- \left[ \intTd{\vr \vu \cdot \bfphi }  \right]_{t=0}^{t=\tau}\\
		&+ \int_0^\tau \intTd{ \Big( \vr \vu \cdot \partial_t \bfphi + \vr \vu \otimes \vu: \Grad \bfphi + p(\chi, \vr) \Div \bfphi - \mathbb{S} : \Ds \bfphi \Big) } \dt \\ & 
		+ \kappa \int_0^\tau \int_{\Td \times S^{d-1}}  \left< \delta V; \bfphi \right>  \ \D V \dt,
	\end{split}
\]	
where, by virtue of \eqref{ds10}, 
\[
\begin{split}
	\int_0^\tau &\intTd{ \Big( F(\chi, \Ds \bfphi) - F(\chi, \Ds \vu) - \mathbb{S} : (\Ds \bfphi - \Ds \vu \Big) }  \dt\\ 
	&= -\int_0^\tau 
	\intTd{ \Big( F(\chi, \Ds \bfphi) - F(\chi, \Ds \vu) - \mathbb{S} : (\Ds \bfphi - \Ds \vu \Big) }  \dt \leq 0.
\end{split}
\]
Consequently, we deduce the standard weak varifold formulation of the momentum equation 
\[
\begin{split}
	&\left[ \intTd{\vr \vu \cdot \bfphi }  \right]_{t=0}^{t=\tau}\\
	&= \int_0^\tau \intTd{ \Big( \vr \vu \cdot \partial_t \bfphi + \vr \vu \otimes \vu: \Grad \bfphi + p(\chi, \vr) \Div \bfphi - \mathbb{S} : \Ds \bfphi \Big) } \dt \\ & 
	+ \kappa \int_0^\tau \int_{\Td \times S^{d-1}}  \left< \delta V; \bfphi \right>  \ \D V \dt
\end{split}
\]	
for any $0 \leq \tau \leq T$, and any $\bfphi \in C^1([0,T] \times \Td; R^d))$.

Note that Definition \ref{pD1} can be easily adapted to incompressible fluids just by 
dropping the pressure term and considering solenoidal test functions in \eqref{ds7}. 

\section{Hypotheses and main results}
\label{H}

Our goal is to handle the largest possible class of the pressure--density equations of state as well as the dissipative potentials. 
For technical reasons, however, certain singularity of $F_i$ will be required to keep the system out of the vacuum. 
Our principal technical hypotheses concerning the dissipative potentials read: 
\begin{equation} \label{H1}
F_i : R^{d \times d}_{\rm sym} \to [0, \infty]\ \mbox{convex l.s.c},\ 
F_i(0) = 0,\ 0 \in {\rm int} ( {\rm Dom}[F_1] \cap  {\rm Dom}[F_2] );
\end{equation}
\begin{equation} \label{H2}
	F_i (\mathbb{D}) \ageq  \Big| \mathbb{D} - \frac{1}{d} {\rm tr}[\mathbb{D}] \Big|^\alpha - 1 
	\ \mbox{for some}\ \alpha > \frac{2d}{d +2};
\end{equation}
\begin{equation} \label{H3}
	F_i (\mathbb{D}) = \infty \ \mbox{whenever}\ {\rm tr}[\mathbb{D}] > \Ov{d} \geq 0;
	\end{equation}
$i = 1,2$. Hypothesis \eqref{H2} seems natural to guarantee integrability of the convective term by means of the energy bounds, hypothesis \eqref{H3} yields a uniform bound
\begin{equation} \label{H4}
	|\Div \vu | \leq \Ov{d} \ \mbox{a.a. in}\ (0,T) \times \Omega.
\end{equation}
We refer the reader to \cite{FeLiMa} for the modelling aspects of \eqref{H3}. 

In the context of incompressible fluids, 
the hypothesis \eqref{H4} becomes irrelevant, while \eqref{H1}, \eqref{H2} can be rewritten as 
\begin{equation} \label{H1a}
	F_i : R^{d \times d}_{\rm sym,0} \to [0, \infty]\ \mbox{convex l.s.c},\ 
	F_i(0) = 0,\ 0 \in {\rm int} ( {\rm Dom}[F_1] \cap  {\rm Dom}[F_2] );
\end{equation}
\begin{equation} \label{H2a}
	F_i (\mathbb{D}) \ageq  | \mathbb{D}  |^\alpha - 1 
	\ \mbox{for some}\ \alpha > \frac{2d}{d +2},
\end{equation}
$i=1,2$. As $F_i$ may become infinite, the above class is broad enough to accommodate the rheology of \emph{thick fluids}, 
see Barnes \cite{Barn} and Rodrigues \cite{Rodr}.

As the velocity $\vu$ enjoys the Sobolev regularity \eqref{ds3} and \eqref{H4} holds, the theory of DiPerna--Lions \cite{DL} applies to \eqref{ds5}, \eqref{ds6}. In particular, for the initial data 
\begin{equation} \label{H5}
0 < \underline{\vr} \leq \vr_0 \leq \Ov{\vr},
\end{equation}
the equation of continuity \eqref{ds6} admits a unique renormalized solution $\vr$ belonging to the class 
\[
\vr \in C([0,T]; L^1(\Td)),\ \underline{\vr} \exp \left( - t \Ov{d} \right) \leq \vr(t,x) \leq 
\Ov{\vr} \exp \left(  t \Ov{d} \right) \ \mbox{for any}\ t,x.
\]
By the same token, the phase variables $\chi$, 
\[
\chi \in C([0,T]; L^1(\Td)),\ \mbox{ ranges in the set}\ \{0,1 \}
\]
provided the same is true for the initial data.
In addition, $\chi$ and $\vr$ satisfy the renormalized equation
\begin{equation} \label{H6}
	\begin{split}
	&\left[ \intO{ b(\chi, \vr) \varphi (0, \cdot) } \right]_{t = 0}^{t = \tau} \\ &=
	\int_0^\tau \intO{ \left[ b(\chi, \vr) \partial_t \varphi + b(\chi, \vr) \vu \cdot \Grad \varphi  + \left( b(\chi, \vr) - \frac{ \partial b(\chi, \vr)}{\partial \vr} \vr \right) \Div \vu \varphi \right] } \dt
	\end{split}
\end{equation}
for any $0 \leq \tau \leq T$, any $\varphi \in C^1([0,T] \times \Ov{\Omega})$, and any $b \in C^1_{\rm loc}([0,1] \times (0, \infty))$.

We are ready to state our main results. We start with a mixture of two compressible fluids without surface tension. 

\begin{Theorem} [{\bf Compressible fluids without surface tension}] \label{TH1}
Let $d=2,3$.	
Suppose that $\kappa = 0$ and $p_i \in C^1_c(0,\infty)$ are increasing functions of $\vr$ for $i=1,2$. Let $F_i$, $i=1,2$, satisfy \eqref{H1}--\eqref{H3} with $\alpha \geq \frac{11}{5}$ if $d=3$ and $\alpha \geq 2$ if $d=2$.	
In addition, suppose there is a constant $k > 0$ such that 
\begin{equation} \label{H6a}
\frac{1}{k} F_1 (\mathbb{D}) - k \leq F_2(\mathbb{D}) \leq k \Big( F_1 (\mathbb{D}) + 1 \Big) \ 
\mbox{for any}\ \mathbb{D} \in R^{d\times d}_{\rm sym}.	
	\end{equation}
Let the initial data satisfy
\[
\chi(0, \cdot)= \chi_0 = 1_{\Omega}, \ \Omega_0 \subset \Td \ \mbox{a Lipschitz domain}, 
\]
\[
\vr(0, \cdot )= \vr_0,\ 
0 < \underline{\vr} \leq \vr_0 \leq \Ov{\vr},\ \vr \vu (0, \cdot) = \vm_0 \in L^2(\Td; R^d).
\]

Then the problem
\eqref{p0}, \eqref{p1}, \eqref{p1a}, \eqref{p5a} admits a dissipative varifold solution $(\vr, \vu, \chi)$ in $(0,T) \times \Td$ in the sense of Definition \ref{pD1}.

	\end{Theorem}

The hypothesis \eqref{H6a} requires the dissipative potentials to share the growth of the same order for large $\mathbb{D}$.

To handle the compressible case with surface tension, we restrict ourselves to the pressure--density rheological law pertinent to isothermal gases, namely 
\begin{equation} \label{gas}
p_i(\vr) = a_i \vr, \ a_i > 0,\ i=1,2,\ p(\vr, \chi) = \chi a_1 \vr + (1 - \chi) a_2 \vr. 	
	\end{equation}

\begin{Theorem} [{\bf Compressible fluids with surface tension}] \label{TH3}
	Let $d=2,3$, and $\kappa \geq 0$. Suppose that the pressure is given by \eqref{gas}. Let $F_i$, $i=1,2$, satisfy \eqref{H1}--\eqref{H3}.
	Let the initial data satisfy
	\[
	\chi(0, \cdot)= \chi_0 = 1_{\Omega}, \ \Omega_0 \subset \Td \ \mbox{a Lipschitz domain}, 
	\]
	\[
	\vr(0, \cdot )= \vr_0,\ 
	0 < \underline{\vr} \leq \vr_0 \leq \Ov{\vr},\ \vr \vu (0, \cdot) = \vm_0 \in L^2(\Td; R^d).
	\]
	
	Then the problem
	\eqref{p0}, \eqref{p1}, \eqref{p1a}, \eqref{p5a} admits a dissipative varifold solution $(\vr, \vu, \chi,V)$ in $(0,T) \times \Td$ in the sense of Definition \ref{pD1}.
	
\end{Theorem}

Finally, we reformulate the result in terms of incompressible fluids.
 
\begin{Theorem} [{\bf Incompressible fluids}] \label{TH2}

Let $d=2,3$, and $\kappa \geq 0$. Suppose that $F_i$ satisfy \eqref{H1a}, \eqref{H2a}. 	
Let the initial data satisfy
\[
\chi(0, \cdot)= \chi_0 = 1_{\Omega}, \ \Omega_0 \subset \Td \ \mbox{a Lipschitz domain}, 
\]
\[
\vr(0, \cdot )= \vr_0,\ 
0 < \underline{\vr} \leq \vr_0 \leq \Ov{\vr},\ \vr \vu (0, \cdot) = \vm_0 \in L^2(\Td; R^d).
\]

Then the problem
\eqref{p0}, \eqref{p1}, \eqref{p1a}, \eqref{p5a} admits a dissipative varifold solution $(\vr, \vu, \chi, V)$ in $(0,T) \times \Td$ in the incompressible setting of Definition \ref{pD1}. Specifically, $\Div \vu = 0$ and \eqref{ds7} holds 
for any test function $\bfphi$ such that $\Div \bfphi = 0$.

\end{Theorem}

The rest of the paper is devoted to the proof of the above results.

\section{Basic approximation scheme}
\label{b}

To construct the varifold solution, we adapt the approximation scheme introduced in \cite{AbbFeiNov}. First, we replace 
the dissipative potentials $F_i$ by their Moreau--Yosida approximation, 
\[
F^\ep_i (\mathbb{D}) = \inf_{\mathbb{M} \in R^{d \times d}_{\rm sym}} \left\{ \frac{1}{2 \ep} | \mathbb{M} - \mathbb{D} |^2 + F^i(\mathbb{D})\right\},\ i=1,2, \ \delta > 0
\]

Next, we consider an orthogonal basis $\{ \vc{e}_n \}_{n=1}^\infty$ of the space $L^2(\Omega; R^d)$ consisting of, say, 
trigonometric polynomials. We look for approximate velocity field belonging to the space 
\[
\vu \in C([0,T]; X_N),\ X_N = {\rm span} \{ \vc{e}_n \}_{n=1}^N .
\]

Given $\vu \in C([0,T]; X_N)$, the transport equation as well as the equation of continuity may be solved by the method of characteristics, 
\begin{equation} \label{b1}
\chi(t, x) = \chi_0 (\vc{X}^{-1}(t,x)),\ 
\vr(t,x) = \vr_0 (\vc{X}^{-1}(t,x)) \exp \left( - \int_0^t \Div \vc{u} (s, \vc{X}^{-1}(s, x )) \D s \right),
\end{equation}
where $\vc{X}$ is the associated Lagrangian flow 
\[
\frac{\D}{\dt} \vc{X}(t,x) = \vu (t, \vc{X}(t,x)),\ \vc{X}(0,x) = x.
\]

At this stage, we approximate the initial data $\chi_0$ to be a characteristic function of a $C^1-$domain $\Omega_0$.
Accordingly, if $\vu$ is smooth with respect to the $x-$variable, the image of $\Omega_0$ under the flow remains $C^1$ at any time, 
\[
\chi_0 = 1_{\Omega_0} \ \mbox{of class}\ C^1 \ \Rightarrow  \ \chi(t, \cdot) = 1_{\Omega_t},\ \Omega_t \ \mbox{of class}\ C^1. 
\]

The velocity field will be identified via a Faedo--Galerkin approximation:
\begin{equation} \label{b2}
	\begin{split}
	\left[ \intTd{ \vr \vu \cdot \bfphi} \right]_{t=0}^{t=\tau} &= \int_0^\tau \intTd{ \Big( 
		\vr \vu \cdot \partial_t \bfphi + \vr \vu \otimes \vu : \Grad \bfphi + p(\chi, \vr) \Div \bfphi \Big) } \dt \\
	& - \int_0^\tau \intTd{ \partial F^\ep (\chi, \Ds \vu) : \Ds \bfphi } \dt + 
	\kappa \int_0^\tau \int_{\partial \Omega_t} \left( \mathbb{I} - \vc{n} \otimes \vc{n} \right) : \Grad \bfphi \ \D S_x \dt\\
	&- \delta \int_0^\tau \intTd{ \Del^m \vu \cdot \Del^m \bfphi } \dt,\ \vr \vu(0, \cdot) = \vm_0,\ \ep > 0,\ m > 2d, 
\end{split}
\end{equation}
for any $\bfphi \in C^1([0,T]; X_N)$,
where we have set 
\[
\partial F^\ep (\Ds \vu) = \chi \partial F^\ep_1 (\Ds \vu) + (1 - \chi) \partial F^\ep_2 (\Ds \vu). 
\]
The reader may consult Abels \cite{AbelsI}, where the specific form of the ``varifold term'' in \eqref{b2} is discussed. 
Identity \eqref{b2} contains also an elliptic regularization necessary for performing the passage from smooth to 
rough dissipation potentials.

The approximation scheme depends on three parameters: $N$, $\ep$, and $\delta$. These being fixed, the existence of 
approximate solutions can be shown by the standard fixed--point argument, see \cite[Section 3]{AbbFeiNov} and Abels \cite[Section 4]{AbelsI}. Our goal is to perform consecutively the limits $N \to \infty$, $\ep \to 0$, and $\delta \to 0$.

\subsection{Energy estimates}

As all quantities are smooth at the basic approximation level, we may consider $\vu$ as a test function in \eqref{b2} obtaining the energy balance
\begin{equation} \label{b3}
	\begin{split}
	\frac{\D}{\dt} \left[ \intTd{ \left( \frac{1}{2} \vr |\vu|^2 + P(\chi,\vr) \right) } 
		+ \kappa \| \Grad \chi \|_{\mathcal{M}^+(\Td)} \right] \\ 
		+ \intTd{ \partial F^\ep (\chi, \Ds \vu) : \Ds \vu } + \delta \intTd{ |\Del^m \vu |^2 } = 0,
		\end{split}
		\end{equation}
where we have used the identity (see Abels \cite[Section 2.4, formula (2.8)]{AbelsI}) 
\begin{equation} \label{b4}
	\frac{\D}{\dt} \int_{\partial \Omega_t} \varphi \ \D S_x = 
	\int_{\partial \Omega_t} (\mathbb{I} - \vc{n} \otimes \vc{n}) : \Grad (\varphi \vu) \ \D S_x
	+ \int_{\partial \Omega_t} \vc{n} \cdot \ \Grad \varphi  \vc{n} \cdot \vu \ \D S_x, 
	\ \varphi \in C^1(\Td),
	\end{equation}
together with the renormalized version \eqref{H6} of the equations \eqref{ds5}, \eqref{ds6}, specifically 
\[
	\left[ \intTd{ P(\chi, \vr) } \right]_{t=0}^{t = \tau} + \int_0^\tau \intTd{p(\chi, \vr) \Div \vu } \dt = 0
	\]
for any $0 \leq \tau \leq T$.

\subsection{Limit $N \to \infty$}
\label{N}

The limit $N \to \infty$ in the family of Faedo--Galerkin approximations is straightforward. 
Let $(\chi_N, \vr_N, \vu_N)_{N \geq 1}$ be the corresponding family of solutions. Keeping 
$\ep > 0$, $\delta > 0$ fixed we deduce from the energy balance and \eqref{b1} the following uniform bounds:
\begin{equation} \label{b6}
	\begin{split}
		\vr_N &\approx 1,\\ 
		(\vu_N)_{N \geq 1} \ &\mbox{bounded in}\ L^\infty (0,T; L^2(\Td; R^d)) 
		\cap L^2(0,T; W^{2m,2} (\Td; R^d )), \\
		(\chi_N)_{N \geq 1}\ &\mbox{bounded in}\ L^\infty((0,T) \times \Td),\  
		(\kappa \chi_N)_{N \geq 1} \ \mbox{bounded in}\ L^\infty(0,T; BV(\Td)).
		\end{split}
\end{equation}

First, we have 
\begin{equation} \label{b7}
	\vu_N \to \vu \ \mbox{weakly-(*) in} \ L^\infty(0,T; L^2(\Td; R^d)) 
	\ \mbox{and weakly in}\ L^2(0,T; W^{2m, 2} (\Td; R^d))
\end{equation}
passing to a suitable subsequence as the case may be. In particular, the velocity 
field remains regular at least in the spatial variable, and we get 
\begin{equation} \label{b8}
\vr_N \to \vr \ \mbox{in} \ C^1([0,T] \times \Td),\ \chi_N \to \chi \ \mbox{weakly-(*) in} \ 
 L^\infty_{{\rm weak-(*)}} (0,T; \mathcal{M}(R^d)),
\end{equation}
and 
\[
\chi_N \to \chi \ \mbox{in}\ C([0,T]; L^1(\Td)) ,
\]
where the limit functions are given by formula \eqref{b1}.

Denoting 
\[
\Pi_N : L^2(\Td; R^d)  \to X_N 
\]
the orthogonal projection, we deduce from \eqref{b2} that 
\[
( \partial_t \Pi_N (\vr_N \vu_N)  )_{N \geq 1} 
\ \mbox{bounded in}\ L^p(0,T; W^{-k,2}(\Td; R^d) ) \ \mbox{for some}\ p > 1,\ k > 0. 
\]
Consequently, by virtue of Aubin--Lions Lemma,
\[
\int_0^T \intTd{ \vr_N |\vu_N|^2 } \dt \to \int_0^T \intTd{ \vr |\vu|^2} \dt,
\]
yielding 
\[
\vu_N \to \vu \ \mbox{in}\ L^2(0,T; L^2(\Td; R^d)).
\]
Thus, finally, 
\[
\vu_N \to \vu \ \mbox{in}\ L^q(0,T; C^1(\Td; R^d)) \ \mbox{for some}\ q > 2.
\]
As the regularization $\partial F^\ep_i$ are globally Lipschitz in $\Ds \vu$, we are allowed to perform the limit in the momentum equation \eqref{b2} obtaining 
\begin{equation} \label{b9}
	\begin{split}
		\left[ \intTd{ \vr \vu \cdot \bfphi} \right]_{t=0}^{t=\tau} &= \int_0^\tau \intTd{ \Big( 
			\vr \vu \cdot \partial_t \bfphi + \vr \vu \otimes \vu : \Grad \bfphi + p(\chi, \vr) \Div \bfphi \Big) } \dt \\
		& - \int_0^\tau \intTd{ \partial F^\ep (\chi, \Ds \vu) : \Ds \bfphi } \dt + 
		\kappa \int_0^\tau \int_{\partial \Omega_t} \left( \mathbb{I} - \vc{n} \otimes \vc{n} \right) : \Grad \bfphi \ \D S_x \dt\\
		&- \delta \int_0^\tau \intTd{ \Del^m \vu \cdot \Del^m \bfphi } \dt,\ \vr \vu(0, \cdot) = \vm_0,\ \ep > 0,\ m > 2d, 
	\end{split}
\end{equation}
for any $\bfphi \in C^1([0,T] \times \Td; R^d)$.

Finally, we pass to the limit in the energy balance \eqref{b3}: 
\begin{equation} \label{b10}
	\begin{split}
		\left[ \psi \intTd{ \left( \frac{1}{2} \vr |\vu|^2 + P(\chi,\vr) \right) } 
		+ \kappa \| \Grad \chi \|_{\mathcal{M}^+(\Td)} \right]_{t = 0}^{ t = \tau} \\ 
		- \int_0^\tau \partial_t \psi \left( \intTd{ \left( \frac{1}{2} \vr |\vu|^2 + P(\chi,\vr) \right) } + \kappa \| \Grad \chi \|_{\mathcal{M}^+(\Td)} \right) \dt\\
		+ \int_0^\tau \psi \intTd{ \partial F^\ep (\chi, \Ds \vu) : \Ds \vu } \dt + \delta \int_0^\tau \psi \intTd{ |\Del^m \vu |^2 } \dt \leq 0,
	\end{split}
\end{equation}
 for a.a. $0 \leq \tau \leq T$ and any $\psi \in C^1[0,T]$, $\psi \geq 0$.
 
 \begin{Remark} \label{bR1}
 	As a matter of fact, a more refined argument would yield equality in \eqref{b10}. This is, however, not needed for the remaining part of the proof.
 	
 	\end{Remark}
 
\section{Compressible case without surface tension}
\label{C} 

Being given a family of approximate solutions identified in Section \ref{N}, our goal is to perform the limits $\ep \to 0$, $\delta \to 0$. Fixing $\delta > 0$
we denote $(\vr_\ep, \chi_\ep, \vu_\ep)_{\ep > 0}$ the associated approximate solution.

\subsection{Limit $\ep \to 0$}
\label{SE}

Since the energy balance \eqref{b10} holds, we may repeat step by step the arguments of Section \ref{N} to obtain 
\begin{equation} \label{C1}
	\begin{split}
		\vr_\ep \approx 1,\ \vre &\to \vr \ \mbox{in}\ C^1([0,T] \times \Omega),\\ 
		\vu_\ep &\to \vu \ \mbox{weakly-(*) in} \ L^\infty(0,T; L^2(\Td; R^d)) 
		\ \mbox{and weakly in}\ L^2(0,T; W^{2m, 2} (\Td; R^d)), \\ 
		\chi_\ep &\to \chi \ \mbox{in}\ C([0,T]; L^1(\Td)) 
		\ \mbox{and weakly-(*) in}\ L^\infty(0,T; BV(\Td))
	\end{split}
	\end{equation}
passing to a suitable subsequence as the case may be. Similarly to the preceding section, the characteristic curves are still well defined and the limit functions $\chi$, $\vr$ are given by formula \eqref{b1}. 

Next, we take $\psi = 1$ in the energy inequality \eqref{b10} and subtract the resulting expression from the momentum balance obtaining
\begin{equation} \label{C2}
	\begin{split}
		&\left[ \intTd{ \left( \frac{1}{2} \vre |\vue|^2 - \vre \vue \cdot \bfphi  + P(\chi_\ep,\vre) \right) } 
		+ \kappa \| \Grad \chi_\ep \|_{\mathcal{M}^+(\Td)} \right]_{t = 0}^{ t = \tau} \\ 
		&+ \int_0^\tau \intTd{ \Big( 
			\vre \vue \cdot \partial_t \bfphi + \vre \vue \otimes \vue : \Grad \bfphi + p(\chi_\ep, \vre) \Div \bfphi \Big) } \dt\\
		&+ \kappa \int_0^\tau \int_{\partial \Omega_{\ep,t}} \left( \mathbb{I} - \vc{n} \otimes \vc{n} \right) : \Grad \bfphi \ \D S_x \dt  - \delta \int_0^\tau \intTd{ \Del^m \vue \cdot \Del^m \bfphi } \dt
		\\
		&+ \int_0^\tau  \intTd{ \partial F^\ep (\chi_\ep , \Ds \vue) : (\Ds \vue - \Ds \bfphi ) } \dt + \delta \int_0^\tau \psi \intTd{ |\Del^m \vue |^2 } \dt \leq 0,
	\end{split}
\end{equation}
for any $\bfphi \in C^1([0,T] \times \Td; R^d)$,
where 
\[
 \intTd{ \partial F^\ep (\chi_\ep , \Ds \vue) : (\Ds \vue - \Ds \bfphi ) }  \geq \intTd{ \Big( F^\ep (\chi_\ep, \Ds \vue) - F^\ep (\chi_\ep, \Ds \bfphi) \Big) }
\]

Now, consider $\phi \in C^1(\Td; R^d)$ such that 
$|\Ds \phi | < M$,
and $\psi \in C^1_c(0,T)$, 
$0 \leq \psi \leq 1$.
In view if hypothesis \eqref{H1}, 
\[
\Ds \phi \in {\rm int} ({\rm Dom}[F_1] \cap {\rm Dom}[F_2])
\]
as long as $M > 0$ is small enough.
Using $\bfphi(t,x) = \psi(t) \phi (x)$ as a test function in \eqref{C2} we obtain 
\[
	\begin{split}
	&\int_0^T \partial_t \psi \intTd{ 
			\vre \vue \cdot \phi } \dt + \int_0^T \psi \intTd{ \Big( \vre \vue \otimes \vue : \Grad \phi + p(\chi_\ep, \vre) \Div \phi \Big) } \dt\\
		&+ \kappa \int_0^T \psi \int_{\partial \Omega_{\ep,t}} \left( \mathbb{I} - \vc{n} \otimes \vc{n} \right) : \Grad \phi \ \D S_x \dt  - \delta \int_0^T \psi \intTd{ \Del^m \vue \cdot \Del^m \phi } \dt
		\\
		 &\leq  \intTd{ \left( \frac{1}{2} \vr_0 |\vu_0|^2  + P(\chi_0,\vr_0) \right) } 
		+ \kappa \| \Grad \chi_0 \|_{\mathcal{M}^+(\Td)} + \int_0^T \intTd{ F^\ep (\chi_\ep, \psi \Ds \phi ) },
	\end{split}
\]
where, furthermore, 
\[
\intTd{ F^\ep (\chi_\ep, \psi \Ds \phi )} \leq \intTd{ F_1 (\psi \Ds \phi) } + 
\intTd{ F_2 (\psi \Ds \phi)  } \leq C(M).
\]
Thus we may infer that 
\[
\begin{split}
t &\in [0,T] \mapsto \intTd{(\vre \vue) (t, \cdot) \cdot \phi } \in BV[0,T] \\
& \mbox{whenever}\ \phi \in C^1(\Td; R^d),
\end{split}
\] 
with the norm bounded uniformly for $\ep \to 0$. Consequently, similarly to Section \ref{N}, we conclude 
\[
\int_0^T \intTd{ \vre |\vue|^2 } \dt \to \int_0^T \intTd{ \vr |\vu|^2 },
\]
yielding, finally,
\begin{equation} \label{C3}
\vue \to \vu \ \mbox{in}\ L^q(0,T; C^1(\Td; R^d)) \ \mbox{for some}\ q > 2.
\end{equation}

The ultimate goal of this step is to perform the limit in \eqref{C2}. Summing up the previous observations, we get 
\begin{equation} \label{C5}
	\begin{split}
		\limsup_{\ep \to 0} &\int_0^\tau \intTd{ F^\ep (\chi_\ep, \Ds \bfphi) } \dt - 
		\liminf_{\ep \to 0}
		\int_0^\tau \intTd{ F^\ep (\chi_\ep, \Ds \vue) } \dt \\
		&\geq \left[ \intTd{ \left( \frac{1}{2} \vr |\vu|^2 - \vr \vu \cdot \bfphi  + P(\chi,\vr) \right) } 
		+ \kappa \| \Grad \chi \|_{\mathcal{M}^+(\Td)} \right]_{t = 0}^{ t = \tau} \\ 
		&+ \int_0^\tau \intTd{ \Big( 
			\vr \vu \cdot \partial_t \bfphi + \vr \vu \otimes \vu : \Grad \bfphi + p(\chi, \vr) \Div \bfphi \Big) } \dt\\
		&+ \kappa \int_0^\tau \int_{\partial \Omega_{t}} \left( \mathbb{I} - \vc{n} \otimes \vc{n} \right) : \Grad \bfphi \ \D S_x \dt  - \delta \int_0^\tau \intTd{ \Del^m \vu \cdot \Del^m \bfphi } \dt
		\\
		&+ \delta \int_0^\tau \psi \intTd{ |\Del^m \vu |^2 } \dt
	\end{split}
\end{equation}
for any $\bfphi \in C^1([0,T] \times \Td; R^d)$. Consider a test function 
\begin{equation} \label{C4}
\bfphi \in C^1([0,T] \times \Td; R^d),\ \int_0^T \intTd{ \Big( F_1(\Ds \bfphi) + F_2 (\Ds \bfphi) \Big) } \dt < \infty,
\end{equation}
and write
\[
\int_0^\tau \intTd{ F^\ep (\chi_\ep, \Ds \bfphi) } \dt = 
\int_0^\tau \intTd{ \chi_\ep F^\ep_1 (\Ds \bfphi) } \dt + 
\int_0^\tau \intTd{ (1 - \chi_\ep) F^\ep_2 (\Ds \bfphi) } \dt.
\]
We already know that
\[
\chi_\ep \to \chi \in \{ 0 ; 1 \}\ \mbox{a.a. in}\ (0,T) \times \Td,\ 
F^\ep_i (\Ds \bfphi) \nearrow F_i(\bfphi) \ \mbox{in}\ (0,T) \times \Td,\ i=1,2.
\]
Consequently, in accordance with \eqref{C4},
\begin{equation} \label{C5a}
\begin{split}
\int_0^\tau &\intTd{ F^\ep (\chi_\ep, \Ds \bfphi) } \dt = 
\int_0^\tau \intTd{ \chi_\ep F^\ep_1 (\Ds \bfphi) } \dt + 
\int_0^\tau \intTd{ (1 - \chi_\ep) F^\ep_2 (\Ds \bfphi) } \dt\\
&\leq \int_0^\tau \intTd{ \chi_\ep F_1 (\Ds \bfphi) } \dt + 
\int_0^\tau \intTd{ (1 - \chi_\ep) F_2 (\Ds \bfphi) } \dt\\
&\to \int_0^\tau \intTd{ \chi F_1 (\Ds \bfphi) } \dt + 
\int_0^\tau \intTd{ (1 - \chi) F_2 (\Ds \bfphi) } \dt = \int_0^\tau \intTd{ F(\chi, \Ds \bfphi) } \dt.
\end{split}
\end{equation}

On the other hand,
\[
\begin{split}
\int_0^\tau &\intTd{F^\ep (\chi_\ep, \Ds \vue) } \dt = 
\int_0^\tau \intTd{ \chi_\ep F^{\ep}_1 (\Ds \vue) } \dt + \int_0^\tau \intTd{ (1 - \chi_\ep) F^{\ep}_2 (\Ds \vue) } \dt\\
&\geq \int_0^\tau \intTd{ \chi_\ep F^{\Ov{\ep}}_1 (\Ds \vue) } \dt + \int_0^\tau \intTd{ (1 - \chi_\ep) F^{\Ov{\ep}}_2 (\Ds \vue) } \dt\\
&\to \int_0^\tau \intTd{ \chi F^{\Ov{\ep}}_1 (\Ds \vu) } \dt + \int_0^\tau \intTd{ (1 - \chi) F^{\Ov{\ep}}_2 (\Ds \vu) } \dt \ \mbox{as}\ \ep \to, \ \Ov{\ep} \ \mbox{fixed}.
\end{split}
\]
Consequently, 
\[
\begin{split}
\liminf_{\ep \to 0} &\int_0^\tau \intTd{F^\ep (\chi_\ep, \Ds \vue) } \dt \\ & \geq
\int_0^\tau \intTd{ \chi F^{\Ov{\ep}}_1 (\Ds \vu) } \dt + \int_0^\tau \intTd{ (1 - \chi) F^{\Ov{\ep}}_2 (\Ds \vu) } \dt
\end{split}
\] 
for any $\Ov{\ep}$. Thus we may infer that 
\begin{equation} \label{C6}
\begin{split}
	\liminf_{\ep \to 0} &\int_0^\tau \intTd{F^\ep (\chi_\ep, \Ds \vue) } \dt \\ & \geq
	\int_0^\tau \intTd{ \chi F_1 (\Ds \vu) } \dt + \int_0^\tau \intTd{ (1 - \chi) F_2 (\Ds \vu) } \dt.
\end{split}
\end{equation}
Combining \eqref{C5}, \eqref{C5a}, and \eqref{C6}, we obtain 
\begin{equation} \label{C8}
	\begin{split}
			\int_0^\tau &\intTd{ F (\chi, \Ds \bfphi) } \dt - 
		\int_0^\tau \intTd{ F (\chi, \Ds \vu) } \dt\\ 
		&\geq \left[ \intTd{ \left( \frac{1}{2} \vr |\vu|^2 - \vr \vu \cdot \bfphi  + P(\chi,\vr) \right) } 
		+ \kappa \| \Grad \chi \|_{\mathcal{M}^+(\Td)} \right]_{t = 0}^{ t = \tau} \\ 
		&+ \int_0^\tau \intTd{ \Big( 
			\vr \vu \cdot \partial_t \bfphi + \vr \vu \otimes \vu : \Grad \bfphi + p(\chi, \vr) \Div \bfphi \Big) } \dt\\
		&+ \kappa \int_0^\tau \int_{\partial \Omega_{t}} \left( \mathbb{I} - \vc{n} \otimes \vc{n} \right) : \Grad \bfphi \ \D S_x \dt  - \delta \int_0^\tau \intTd{ \Del^m \vu \cdot \Del^m \bfphi } \dt
		\\
		&+ \delta \int_0^\tau \psi \intTd{ |\Del^m \vu |^2 } \dt
	\end{split}
\end{equation}
for any $\bfphi \in C^1([0,T] \times \Td; R^d)$ satisfying \eqref{C4}. 

\subsection{Limit $\delta \to 0$}
\label{SD}

Let $(\chi_\delta, \vrd, \vud )_{\delta > 0}$ be the family of approximate solutions obtained in the previous section. 
Our ultimate goal is to perform the limit $\delta \to 0$. Up to this moment, the proof has been identical for 
Theorems \ref{TH1}, \ref{TH3}, with the obvious modifications for Theorem \ref{TH2}. Now, we focus on the 
case of compressible fluids without surface tension.

If $\kappa = 0$, the approximate solutions satisfy the following  relations:
\begin{equation} \label{C9}
	\left[ \intTd{ \chi_\delta \varphi } \right]_{t=0}^{t = \tau} = 
	\int_0^\tau \intTd{ \Big( \chi_\delta \partial_t \varphi + \chi_\delta \vud \cdot \Grad \varphi + \chi_\delta \Div \vud \varphi \Big) } \dt
\end{equation}
for any $0 \leq \tau \leq T$, $\varphi \in C^1([0,T \times \Td])$; 

\begin{equation} \label{C10}
	\left[ \intTd{ \vrd \varphi } \right]_{t=0}^{t = \tau} = 
	\int_0^\tau \intTd{ \Big( \vrd \partial_t \varphi + \vrd \vud \cdot \Grad \varphi \Big) } \dt
\end{equation}
for any $0 \leq \tau \leq T$, $\varphi \in C^1([0,T] \times \Td)$;
\begin{equation} \label{C11}
	\begin{split}
		\int_0^\tau &\intTd{ F (\chid, \Ds \bfphi) } \dt - 
		\int_0^\tau \intTd{ F (\chid, \Ds \vud) } \dt\\
		&\geq \left[ \intTd{ \left( \frac{1}{2} \vrd |\vud|^2 - \vrd \vud \cdot \bfphi  + P(\chid,\vrd) \right) } 
	 \right]_{t = 0}^{ t = \tau} \\ 
		&+ \int_0^\tau \intTd{ \Big( 
			\vrd \vud \cdot \partial_t \bfphi + \vrd \vud \otimes \vud : \Grad \bfphi + p(\chid, \vrd) \Div \bfphi \Big) } \dt\\
		&- \delta \int_0^\tau \intTd{ \Del^m \vud \cdot \Del^m \bfphi } \dt
		+ \delta \int_0^\tau  \intTd{ |\Del^m \vud |^2 } \dt
	\end{split}
\end{equation}
for any $\bfphi \in C^1([0,T] \times \Td; R^d)$ satisfying \eqref{C4}. Setting $\bfphi \equiv 1$ in \eqref{C11}, we recover the energy estimates
\[
\begin{split}
{\rm ess} \sup_{t \in (0,T)} \left\| \vrd |\vud|^2 \right\|_{L^1(\Td)} &\aleq 1,\\ 
\int_0^T \intTd{ F(\chid, \Ds \vud) } &\aleq 1,\\ 
\delta \int_0^T \intTd | \Del^m \vud |^2 &\aleq 1.
\end{split}
\]
Consequently, by virtue of hypothesis \eqref{H3}, 
\begin{equation} \label{C12a}
| \Div \vud (t,x)| \leq \Ov{d} \ \mbox{for a.a.}\ (t,x).
\end{equation}
In addition, 
\[
\chid \in C([0,T]; L^1(\Td)),\ \chid(t,x) \in \{ 0,1 \} \ \mbox{for a.a.} \ (t,x)
\]
and
\[
\vrd \in C([0,T]; L^1(\Td)),\ \underline{\vr} \exp \left( - t \Ov{d} \right) \leq \vrd(t,x) \leq 
\Ov{\vr} \exp \left(  t \Ov{d} \right) \ \mbox{for a.a}\ (t,x), 
\]
represent the unique renormalized solution of \eqref{C9}, \eqref{C10} in the sense of DiPerna and Lions \cite{DL}.

Similarly to Section \ref{SE}, we show 
\[
\vrd \to \vr \ \mbox{in}\ C_{\rm weak}([0,T]; L^q(\Td)),\ 1 < q < \infty, 
\]
\[
\vrd \vud \to \vr \vu \ \mbox{in} \ L^p_{\rm weak}(0,T; L^2(\Td; R^d)) \ \mbox{for any}\ 1 < p < \infty, 
\]
and 
\begin{equation} \label{C12}
\vrd \vud \otimes \vud \to \vr \vu \otimes \vu \ \mbox{weakly in}\ L^q((0,T) \times \Td; R^{d \times d}_{\rm sym})
\ \mbox{for some}\ q > 1.
\end{equation}
In particular, 
\[
\int_0^T \intTd{\vrd | \vud - \vu |^2  } = \int_0^T \intTd{ \left( \vrd|\vud|^2 - 2 \vrd \vud \cdot \vu + \vrd |\vud|^2 \right) } \dt \to 0. 
\]
As $\vrd$ is bounded below away from zero, the above relation yields 
\begin{equation} \label{C13}
\vud \to \vu\ \mbox{in}\ L^2((0,T) \times \Td; R^d). 
\end{equation}

With the estimates \eqref{C12a}, \eqref{C13} at hand, we may apply \cite[Theorem II.5]{DL}
to obtain 
\begin{equation} \label{C15}
\chid \to \chi \ \mbox{in}\ C([0,T]; L^q(\Td)) \ \mbox{for any}\ 1 \leq q < \infty, 
\end{equation}
where $\chi$ is a renormalized solution of the transport equation \eqref{ds5}, 
$\chi(t,x) \in \{ 0, 1 \}$ for a.a. $(t,x)$. Moreover, it is easy to check that $\vr$ is a renormalized solution of the equation of continuity \eqref{ds6}. 

Finally, we use the fact that $\chid$, $\vrd$ satisfy the renormalized equation \eqref{H6} and rewrite \eqref{C11} in the form  
\begin{equation} \label{C16}
	\begin{split}
		\int_0^\tau &\intTd{ F (\chid, \Ds \bfphi) } \dt - 
		\int_0^\tau \intTd{ F (\chid, \Ds \vud) } \dt\\
		&\geq \left[ \intTd{ \left( \frac{1}{2} \vrd |\vud|^2 - \vrd \vud \cdot \bfphi  \right) } 
		\right]_{t = 0}^{ t = \tau} \\ 
		&+ \int_0^\tau \intTd{ \Big( 
			\vrd \vud \cdot \partial_t \bfphi + \vrd \vud \otimes \vud : \Grad \bfphi + p(\chid, \vrd) (\Div \bfphi - \Div \vud) \Big) } \dt\\
		&- \delta \int_0^\tau \intTd{ \Del^m \vud \cdot \Del^m \bfphi } \dt
		\\
		&+ \delta \int_0^\tau  \intTd{ |\Del^m \vud |^2 } \dt.
	\end{split}
\end{equation}

Repeating the arguments of Section \ref{SE}, we let $\delta \to 0$ in \eqref{C16} obtaining
\begin{equation} \label{C17}
	\begin{split}
			\int_0^\tau &\intTd{ F (\chi, \Ds \bfphi) } \dt - 
		\int_0^\tau \intTd{ F (\chi, \Ds \vu) } \dt \\
		&\geq \left[ \intTd{ \left( \frac{1}{2} \vr |\vu|^2 - \vr \vu \cdot \bfphi  \right) } 
		\right]_{t = 0}^{ t = \tau} \\ 
		&+ \int_0^\tau \intTd{ \Big( 
			\vr \vu \cdot \partial_t \bfphi + \vr \vu \otimes \vu : \Grad \bfphi + \Ov{p(\chi, \vr)} \Div \bfphi - \Ov{p(\chi, \vr) \Div \vu } \Big) } \dt,    
	\end{split}
\end{equation}
where 
\[
\begin{split}
p(\chid, \vrd) &\to \Ov{p(\vr, \chi)} \ \mbox{weakly-(*) in}\ L^\infty((0,T) \times \Td), \\
p(\chid, \vrd) \Div \vud &\to \Ov{p(\vr, \chi) \Div \vud} \ \mbox{weakly in} 
L^\alpha ((0,T) \times \Omega; R^d). 
\end{split}
\]

\subsubsection{Strong convergence of the density}
\label{SCD}

Our ultimate goal in the proof of Theorem \ref{TH1} is to ``remove'' the bars in \eqref{C17} which amounts to showing strong a.a. pointwise convergence 
\begin{equation} \label{C18}
\vrd \to \vr \ \mbox{a.a. in}\ (0,T) \times \Td.
\end{equation}

To see \eqref{C18}, it is enough to justify the choice $\bfphi = \vu$ as a test function 
in \eqref{C17}. Indeed, as shown in \cite{FeLiMa}, this would yield  
\[
\int_0^\tau \intTd{ \Ov{p(\chi, \vr) \Div \vu }
} \dt \geq \int_0^\tau \intTd{ \Ov{p(\chi, \vr)} \Div \vu 
} \dt.
\]
On the other hand, from the renormalization,
\begin{equation} \label{C19}
	\begin{split}
\frac{\D }{\dt} & \intTd{ \Big[ \Ov{ P(\chi, \vr ) } - P(\chi, \vr) \Big] } \dt 
= \intTd{ p(\chi, \vr) \Div \vu - \Ov{ p(\chi, \vr) \Div \vu } } \dt \\ 
&\leq \intTd{ \Big( p(\chi, \vr) - \Ov{p(\chi, \vr)} \Big) \Div \vu } \dt.
\end{split}
\end{equation}
As $P_i$, $i=1,2$, are strictly convex on the range of $\vr$ and $\chid$ converges strongly, we get 
\[
\begin{split}
&\intTd{ \Ov{ P(\chi, \vr ) } - P(\chi, \vr) } = 
\intTd{ \left( \chi \Big[ \Ov{P_1(\vr)} - P_1(\vr) \Big]  + (1 - \chi)
\Big[ \Ov{P_2(\vr)} - P_2(\vr) \Big] \right) } \\ 
&\ageq \limsup_{\delta \to 0} \| \vrd - \vr \|^2_{L^2(\Td)}. 	
\end{split}
\] 
By the same token,
\[ 
\intTd{ \Big( p(\chi, \vr) - \Ov{p(\chi, \vr)} \Big) \Div \vu } \leq 
\| \Div \vu \|_{L^\infty(\Td)} \limsup_{\delta \to 0} \| \vrd - \vr \|^2_{L^2(\Td)}.
\] 
Thus we may apply Gronwall's lemma to \eqref{C19} to obtain, up to a suitable subsequence, 
\eqref{C18}.

In order to make the preceding step rigorous, we have to justify the choice $\bfphi = 
\vu$ in \eqref{C17}. Similarly to \cite{FeLiMa}, we consider the regularization in time via the Steklov averages. Specifically, let 
\[
\begin{split}
\eta_h(t) &= \frac{1}{h} 1_{[-h,0]}(t),\  \eta_{-h}(t) = \frac{1}{h} 1_{[0,h]}(t),\ h > 0,  \\
\xi_\ep &\in \DC (0, \tau),\ 0 \leq \xi_\ep \leq 1,\ \xi_\ep (t) = 1 \ \mbox{if}\ \ep \leq t \leq \tau - \ep.
\end{split}
\]
A suitable approximation of $\vu$ reads 
\[
[\vu]_{h, \ep} = \xi_\ep \eta_{-h} * \eta_h *( \xi_\ep \vu),\ \ep > 0,\ h > 0, 
\] 
where $*$ stands for the convolution in the time variable $t$. The function $[\vu]_{h, \ep}$ can be used as a test function 
in \eqref{C17}. As shown in \cite[Section 4.3]{FeLiMa}, we have 
\[
		\left[ \intTd{ \left( \frac{1}{2} \vr |\vu|^2 - \vr \vu \cdot [\vu_{h,\ep}]  \right) } 
		\right]_{t = 0}^{ t = \tau} + \int_0^\tau \intTd{ \Big( 
			\vr \vu \cdot \partial_t  [\vu_{h,\ep}] + \vr \vu \otimes \vu : \Grad  [\vu_{h,\ep}] } \dt \to 0
\]
as $\ep, h \to 0$. Note that it is only at this point, where the hypothesis $\alpha \geq \frac{11}{5}$ is needed. Consequently, it is enough to show
\begin{equation} \label{C20}
	\limsup_{\ep, h \to 0} \int_0^\tau \intTd{ F(\chi,\Ds [\vu]_{h,\ep} ) } \dt \leq 
\int_0^\tau \intTd{ F(\chi, \Ds \vu ) } \dt.	
\end{equation}
To see \eqref{C20} we first use convexity of $F_i$ and the hypothesis $F_i(0) = 0$ to observe 
\[
F_i (\Ds [\vu]_{h,\ep} ) ) = F_i (\xi_\ep \eta_{-h}* \eta_h * (\xi_\ep \Ds \vu)) \leq 
\xi_\ep F_i (\eta_{-h}* \eta_h * (\xi_\ep \Ds \vu)),\ i = 1,2.
\]
Next, by Jensen' inequality, 
\[
F_i (\eta_{-h}* \eta_h * (\xi_\ep \Ds \vu)) \leq \eta_{-h} * F_i (\eta_h * (\xi_\ep \Ds \vu)) \leq \eta_{-h} *
\eta_h * F_i (\xi_\ep \Ds \vu)),\ i=1,2.
\]
Consequently, 
\[
F_i (\Ds [\vu]_{h,\ep} ) ) \leq \xi_\ep \eta_{-h} *
\eta_h * F_i (\xi_\ep \Ds \vu)). 
\]
By virtue of hypothesis \eqref{H6a},
\[
F_i (\xi_\ep \Ds \vu) \to F_i (\Ds \vu) \ \mbox{in}\ L^1 ((0,\tau) \times \Td),\ i=1,2,\ \mbox{as}\ \ep \to 0; 
\]
whence \eqref{C20} follows. Note that it is only at this moment, where the hypothesis \eqref{H6a} was needed. 

We have proved Theorem \ref{TH1}.

\section{Fluids with surface tension}

The above arguments can be easily adapted to handle the general case of compressible fluids with surface tension stated 
in Theorem \ref{TH3}. Indeed, as the pressure $p(\chi, \vr)$ is now linear with respect to $\vr$, the proof of strong 
convergence of the density performed in Section \ref{SCD} is now longer needed. Accordingly, we can drop the hypothesis 
$\alpha \geq \frac{11}{5}$ as well as \eqref{H6a}.

In addition to the arguments of the preceding section, we have to perform the limit $\delta \to 0$ in the varifold term 
\[
\kappa \int_0^\tau \int_{\partial \Omega_{\delta, t}} \left( \mathbb{I} - \vc{n} \otimes \vc{n} \right) : \Grad \bfphi \ \D S_x \dt
\]
and to show the compatibility relation \eqref{ds8}. However, this can be done in the same way is in \cite{AbelsI}, \cite{Plot2}.
This completes the proof of Theorem \ref{TH3}

Finally, it is easy to observe that the proof can be modified in a straightforward manner to deal with the incompressible case claimed in Theorem \ref{TH2}.

\def\cprime{$'$} \def\ocirc#1{\ifmmode\setbox0=\hbox{$#1$}\dimen0=\ht0
	\advance\dimen0 by1pt\rlap{\hbox to\wd0{\hss\raise\dimen0
			\hbox{\hskip.2em$\scriptscriptstyle\circ$}\hss}}#1\else {\accent"17 #1}\fi}


\begin{thebibliography}{10}
	
	\bibitem{AbbFeiNov}
	A.~Abbatiello, E.~Feireisl, and A.~Novotn\'{y}.
	\newblock Generalized solutions to models of compressible viscous fluids.
	\newblock {\em Discrete Contin. Dyn. Syst.}, {\bf 41}(1):1--28, 2021.
	
	\bibitem{AbelsI}
	H.~Abels.
	\newblock On generalized solutions of two-phase flows for viscous
	incompressible fluids.
	\newblock {\em Interfaces Free Bound.}, {\ bf 9}(1):31--65, 2007.
	
	\bibitem{AbelsII}
	H.~Abels.
	\newblock On the notion of generalized solutions of viscous incompressible
	two-phase flows.
	\newblock In {\em Kyoto {C}onference on the {N}avier-{S}tokes {E}quations and
		their {A}pplications}, RIMS K\^{o}ky\^{u}roku Bessatsu, B1, pages 1--19. Res.
	Inst. Math. Sci. (RIMS), Kyoto, 2007.
	
	\bibitem{ALNS}
	D.~M. Ambrose, M.~C. Lopes~Filho, H.~J. Nussenzveig~Lopes, and W.~A. Strauss.
	\newblock Transport of interfaces with surface tension by 2{D} viscous flows.
	\newblock {\em Interfaces Free Bound.}, {\bf 12}(1):23--44, 2010.
	
	\bibitem{Barn}
	H.~A. Barnes.
	\newblock Shear-{T}hickening (“{D}ilatancy”) in suspensions on
	nonaggregating solidparticles dispersed in {N}ewtonian liquids.
	\newblock {\em J. Rheology}, {\bf 33}:329--366, 1989.
	
	\bibitem{BuGwMaSG}
	M.~Bul\'{\i}\v{c}ek, P.~Gwiazda, J.~M\'{a}lek, and A.~\'{S}wierczewska Gwiazda.
	\newblock On unsteady flows of implicitly constituted incompressible fluids.
	\newblock {\em SIAM J. Math. Anal.}, {\bf 44}(4):2756--2801, 2012.
	
	\bibitem{DeSo}
I. Denisova, V. A. Solonnikov
\newblock{Local and global solvability of free boundary value problems near equalibria}
\newblock{Handbook of Mathematical Analysis in Mechanics of Viscous Fluids, Vol. 2}
\newblock{eds. Y. Giga, A. Novotny},  Springer, 2018
	
	\bibitem{DL}
	R.J. DiPerna and P.-L. Lions.
	\newblock Ordinary differential equations, transport theory and {S}obolev
	spaces.
	\newblock {\em Invent. Math.}, {\bf 98}:511--547, 1989.
	
	\bibitem{FeLiMa}
	E.~Feireisl, X.~Liao, and J.~M{\' a}lek.
	\newblock Global weak solutions to a class of non-{N}ewtonian compressible
	fluids.
	\newblock {\em Math. Meth. Appl. Sci}, {\bf 38}:3482--3494, 2015.
	
	\bibitem{Plot1}
	P.~I. Plotnikov.
	\newblock Generalized solutions of a problem on the motion of a non-{N}ewtonian
	fluid with a free boundary.
	\newblock {\em Sibirsk. Mat. Zh.}, {\bf 34}(4):127--141, iii, ix, 1993.
	
	\bibitem{Plot3}
	P.~I. Plotnikov.
	\newblock Varifold solutions of a free boundary problem in viscous fluid
	dynamics.
	\newblock In {\em Free boundary problems in fluid flow with applications
		({M}ontreal, {PQ}, 1990)}, volume 282 of {\em Pitman Res. Notes Math. Ser.},
	pages 28--32. Longman Sci. Tech., Harlow, 1993.
	
	\bibitem{Plot2}
	P.~I. Plotnikov.
	\newblock Compressible {S}tokes flow driven by capillarity on a free surface.
	\newblock In {\em Navier-{S}tokes equations and related nonlinear problems
		({P}alanga, 1997)}, pages 217--238. VSP, Utrecht, 1998.
	
	\bibitem{Rodr}
	J.-F. Rodrigues.
	\newblock On the mathematical analysis of thick fluids.
	\newblock {\em J. Math. Sci. (N.Y.)}, 210({\bf 6}):835--848, 2015.
	
\end{thebibliography}
\end{document}